\documentclass[11pt]{article}

\usepackage{latexsym,amsmath,amsthm,amsfonts,amsbsy,amscd}
\newtheorem{lemma}{Lemma}
\newtheorem{corollary}{Corollary}

\renewcommand{\pmod}[1]{\; (\bmod\, #1)}

\begin{document}

\title{A note on the trinomial analogue of Bailey's lemma}

\author{\Large
\Large S.~Ole Warnaar\thanks{
e-mail: {\tt warnaar@maths.mu.oz.au}} \\
\mbox{} \\
\em Department of Mathematics, University of Melbourne\\
\em Parkville, Victoria 3052, Australia}
 
\date{\Large February, 1997}
\maketitle

\begin{abstract}
Recently, Andrews and Berkovich introduced a trinomial version of
Bailey's lemma. In this note we show that each ordinary Bailey pair 
gives rise to a trinomial Bailey pair. This largely widens the
applicability of the trinomial Bailey lemma and proves some of the 
identities proposed by Andrews and Berkovich.
\end{abstract}

\subsection*{The trinomial Bailey lemma}
In a recent paper, Andrews and Berkovich (AB) proposed a trinomial analogue of
Bailey's lemma~\cite{AB}. 
As starting point AB take the following definitions of the $q$-trinomial 
coefficients
\begin{equation}
\binom{L;B;q}{A}_2 = \sum_{j=0}^{\infty} \frac{q^{j(j+B)}(q)_L}
{(q)_j(q)_{j+A}(q)_{L-2j-A}}
\end{equation}
and
\begin{equation}
T_n(L,A,q) = q^{\frac{L(L-n)-A(A-n)}{2}} \binom{L;A-n;q^{-1}}{A}_2.
\end{equation}
Here $(a)_{\infty}=\prod_{n=0}^{\infty}(1-aq^n)$ and
$(a)_n=(a)_{\infty}/(aq^n)_{\infty}$, $n\in \mathbb{Z}$.
To simplify equations it will also be convenient to introduce the
notation
\begin{equation}
Q_n(L,A,q) = T_n(L,A,q)/(q)_L.
\end{equation}
We note that the $q$-trinomial coefficients are non-zero for $-L\leq A \leq
L$ only.

A pair of sequences 
$\tilde\alpha=\{\tilde\alpha_L\}_{L\geq 0}$ and 
$\tilde\beta=\{\tilde\beta_L\}_{L\geq 0}$ is said to form a
trinomial Bailey pair relative to $n$ if
\begin{equation}
\tilde\beta_L=\sum_{r=0}^L 
Q_n(L,r,q)\, \tilde\alpha_r .
\end{equation}
The trinomial analogue of the Bailey lemma is stated as follows~\cite{AB}.
\begin{lemma}
If $(\tilde\alpha,\tilde\beta)$ is a trinomial Bailey pair relative to $0$,
then
\begin{equation}
\sum_{L=0}^{\infty} (-1)_L \, q^{L/2} \tilde\beta_L =
(-1)_{M+1} \sum_{L=0}^{\infty}
\frac{\tilde\alpha_L}{q^{L/2}+q^{-L/2}} \: Q_1(M,L,q).
\end{equation}
Similarly, if $(\tilde\alpha,\tilde\beta)$ is a trinomial Bailey pair 
relative to $1$, then
\begin{multline}
\sum_{L=0}^{\infty} \big(-q^{-1}\big)_L\, q^L \tilde\beta_L  \\
= (-1)_M \sum_{L=0}^{\infty} \tilde\alpha_L \biggl\{ Q_1(M,L,q) 
-\frac{Q_1(M-1,L+1,q)}{1+q^{-L-1}} 
-\frac{Q_1(M-1,L-1,q)}{1+q^{L-1}} \biggr\}.
\end{multline}
\end{lemma}
As a corollary of their lemma, AB obtain the identities
\begin{equation}\label{Cor1}
\frac{1}{2} \sum_{L=0}^{\infty} (-1)_L \, q^{L/2} \tilde\beta_L =
\frac{(-q)^2_{\infty}}{(q)^2_{\infty}} \sum_{L=0}^{\infty}
\frac{\tilde\alpha_L}{q^{L/2}+q^{-L/2}},
\end{equation}
for a trinomial Bailey pair relative to $0$, and
\begin{equation}\label{Cor2}
\frac{1}{2}\sum_{L=0}^{\infty} \big(-q^{-1}\big)_L\, q^L \tilde\beta_L 
= \frac{(-q)^2_{\infty}}{(q)^2_{\infty}} \sum_{L=0}^{\infty}
\tilde\alpha_L 
\biggl\{\frac{1}{1+q^{L+1}} -\frac{1}{1+q^{L-1}}\biggr\},
\end{equation}
for a trinomial Bailey pair relative to $1$.

\subsection*{From binomial to trinomial Bailey pairs}
In ref.~\cite{AB}, the equations \eqref{Cor1} and \eqref{Cor2}
are used to derive several new $q$-series identities. As input AB take
trinomial Bailey pairs obtained from polynomial identities which
on one side involve $q$-trinomial coefficients.
Among these identities is an identity by the author which was stated
in ref.~\cite{W} without proof, and therefore AB conclude ``We have
checked that his conjecture implies'' followed by their equation (3.21),
which is an identity for the characters of the
$N=2$ superconformal models $SM(2p,(p-1)/2)$.

We now point out that equation (3.21) is a simple consequence of 
lemma~\ref{lemma} stated below. 
First we recall the definition of the ordinary (i.e., binomial) Bailey pair.
A pair of sequences $(\alpha,\beta)$ such that
\begin{equation}\label{BP}
\beta_L=\sum_{r=0}^L 
\frac{\alpha_r}{(q)_{L-r}(aq)_{L+r}}
\end{equation}
is said to form a Bailey pair relative to $a$.

\begin{lemma}\label{lemma}
Let $(\alpha,\beta)$ form a Bailey pair relative to $a=q^{\ell}$,
where $\ell$ is a non-negative integer.
For $n=0,1$, the following identity holds:
\begin{equation}\label{eqlemma}
\sum_{\substack{s=0 \\ s \equiv L+\ell \pmod{2}}}^{L-\ell}
\frac{q^{s(s-n)/2}}{(q)_{\ell} (q)_s} \beta_{(L-s-\ell)/2} = 
\sum_{r=0}^{\infty} Q_n(L,2r+\ell,q) \, \alpha_r .
\end{equation}
\end{lemma}
For $\ell>L$ the above of course trivializes to $0=0$.

Before proving lemma~\ref{lemma} we note an immediate consequence.
\begin{corollary}\label{cor}
Let $(\alpha,\beta)$ form a Bailey pair relative to $a=q^{\ell}$
with non-negative integer $\ell$.
Then $(\tilde\alpha,\tilde\beta)$ defined as
\begin{align}
\tilde\alpha_0,\dots,\tilde\alpha_{\ell-1}=0, &\qquad
\tilde\alpha_{2L+\ell}  = \alpha_L, \qquad
\tilde\alpha_{2L+\ell+1}  = 0, \qquad L\geq 0 \notag \\
\tilde\beta_0,\dots,\tilde\beta_{\ell-1}=0, & \qquad
\tilde\beta_{L+\ell} = 
\sum_{\substack{s=0 \\ s \equiv L \pmod{2}}}^L
\frac{q^{s(s-n)/2}}{(q)_{\ell}(q)_s} \beta_{(L-s)/2}, \qquad L\geq 0
\end{align}
forms a trinomial Bailey pair relative to $n=0,1$.
\end{corollary}

\begin{proof}
The proof is trivial once one adopts the representation of the
$q$-trinomial coefficients as given by equations (2.58) and (2.59)
of ref.~\cite{ABa},
\begin{equation}
Q_n(L,A,q) = \frac{T_n(L,A,q)}{(q)_L} =
\sum_{\substack{s=0 \\ s\equiv L+A \pmod{2}}}^{\infty} 
\frac{q^{s(s-n)/2}}
{(q)_{\frac{L-A-s}{2}}(q)_{\frac{L+A-s}{2}}(q)_s}\: , \qquad n=0,1.
\end{equation}

Now take the defining relation \eqref{BP} of a Bailey pair 
with $a=q^{\ell}$ and make the replacement $L\to (L-s-\ell)/2$ where $s$ 
is an integer $0\leq s \leq L-\ell$ such that $s\equiv L+\ell \pmod{2}$.
After multiplication by $q^{s(s-n)/2}/(q)_s$ this becomes
\begin{equation}
\frac{q^{s(s-n)/2}}{(q)_s} \beta_{(L-s-\ell)/2} = 
(q)_{\ell} \sum_{r=0}^{\infty}
\frac{\alpha_r q^{s(s-n)/2}}
{(q)_{\frac{L-s-\ell}{2}-r} (q)_{\frac{L-s+\ell}{2}+r} (q)_s}.
\end{equation}
Summing over $s$ yields equation \eqref{eqlemma}.
\end{proof}

Returning to AB's paper, we note that their equation (3.21) simply follows
from corollary~\ref{cor} and the ``$M(p-1,p)$ Bailey pairs'' which arises 
from the $M(p-1,p)$ polynomial identities proven in refs.~\cite{B,W}. 
Of course, an equivalent
statement is that the ``conjecture of ref.~\cite{W}'', is proven using
lemma~\ref{lemma} and the $M(p-1,p)$ Bailey pairs.
To make this somewhat more explicit we consider the special case $p=3$.
Then the $M(2,3)$ Bailey pairs are nothing but the entries
A(1) and A(2) of Slater's list~\cite{Sl}.
Specifically, A(1) contains the following Bailey pair relative to $1$:
\begin{equation}
\alpha_L = \begin{cases}
q^{6j^2-j}, & L=3j \geq 0 \\
q^{6j^2+j}, & L=3j >0 \\
-q^{6j^2-5j+1}, & L=3j-1 >0 \\
-q^{6j^2+5j+1}, & L=3j+1 >0
\end{cases}  \qquad \text{and} \qquad
\beta_L  = \frac{1}{(q)_{2L}}.
\end{equation}
By application of corollary~\ref{cor} this gives the trinomial
Bailey pair 
\begin{equation}
\tilde\alpha_L = \begin{cases}
q^{6j^2-j}, & L=6j \geq 0\\
q^{6j^2+j}, & L=6j >0 \\
-q^{6j^2-5j+1}, & L=6j-1 >0 \\
-q^{6j^2+5j+1}, & L=6j+1 >0
\end{cases}   \qquad \text{and} \qquad 
\tilde\beta_L = \! 
\sum_{\substack{s=0 \\ s \equiv L \pmod{2}}}^L \!
\frac{q^{s(s-n)/2}}{(q)_s (q)_{L-s}}.
\end{equation}
Likewise, using entry A(2), we get
\begin{equation}
\tilde\alpha_L = \begin{cases}
q^{6j^2-j},& L=6j-1 >0 \\
q^{6j^2+j},& L=6j+1 >0 \\
-q^{6j^2-5j+1}, & L=6j-3 >0 \\
-q^{6j^2+5j+1}, & L=6j+3 >0
\end{cases} \qquad \text{and} \qquad
\tilde\beta_L  = \!
\sum_{\substack{s=0 \\ s \not\equiv L \pmod{2}}}^L \!
\frac{q^{s(s-n)/2}}{(q)_s(q)_{L-s}}.
\end{equation}
Setting $n=0$ and summing up both trinomial Bailey pairs, we arrive
at the trinomial Bailey pair of equations~(3.18) and (3.19) of~\cite{AB}. 
(Unlike the case $p\geq 4$, this trinomial Bailey pair was actually proven by
AB, using theorem~5.1 of ref.~\cite{A}.)
Very similar results can be obtained through application of Slater's
A(3) and A(4), A(5) and A(6), and A(7) and A(8).

\subsection*{Conclusion}
We conclude this note with several remarks.
First, it is of course not true that each trinomial Bailey pair
is a consequence of an ordinary Bailey pair. The pairs given by
equations (3.13) and (3.14) of ref.~\cite{AB} 
being examples of irreducible trinomial Bailey pairs.
Second, if one replaces $Q_n(L,r,q)$ by its $q$-multinomial
analogue~\cite{S,Wb} and takes that as the definition of a
$q$-multinomial Bailey pair, it becomes straightforward to again
construct multinomial Bailey pairs out of ordinary ones.
Finally, it is worthwhile to note that the Bailey flow from the minimal
model $M(p,p+1)$ to the $N=2$ superconformal model $SM(2p,(p-1)/2)$ as
concluded by AB could now be replaced by $M(p-1,p)\to N=2~SM(2p,(p-1)/2)$.
Perhaps better though would be to write $M(p-1,p)\to M(p,p+1)\to
N=2~SM(2p,(p-1)/2)$, where the first arrow indicates the flow induced by
corollary~\ref{cor} and the second arrow the flow induced by
\eqref{Cor1} and \eqref{Cor2}.

\subsection*{Acknowledgements}
I thank Alexander Berkovich for helpful comments.
This work is supported by the Australian Research Council.


\begin{thebibliography}{99}

\bibitem{A}
G.~E.~Andrews,
{\em Euler's ``Exemplum memorabile inductionis fallacis'' and
$q$-trinomial coefficients}
J. Amer. Math. Soc. {\bf 3} (1990), 653--669.

\bibitem{ABa}
G.~E.~Andrews and R.~J.~Baxter,
{\em Lattice gas generalization of the hard hexagon model. III. 
$q$-Trinomial coefficients},
J. Stat. Phys. {\bf 47} (1987), 297--330.

\bibitem{AB}
G.~E.~Andrews and A.~Berkovich,
{\em A trinomial analogue of Bailey's lemma and $N=2$ superconformal
invariance}, q-alg/9702008. Submitted to Commun. Math. Phys.

\bibitem{B}
A.~Berkovich,
{\em Fermionic counting of RSOS-states and Virasoro character formulas
for the unitary minimal series $M(\nu,\nu+1)$. Exact results},
Nucl. Phys. B {\bf 431} (1994), 315--348.

\bibitem{S}
A.~Schilling,
{\em Multinomials and polynomial bosonic forms for the branching functions
of the $\widehat{su}(2)_M \times \widehat{su}(2)_N / \widehat{su}(2)_{N+M}$
conformal coset models},
Nucl. Phys. B {\bf 467} (1996), 247--271.

\bibitem{Sl}
L.~J.~Slater,
{\em A new proof of Rogers's transformations of infinite series},
Proc. London Math. Soc. (2) {\bf 53} (1951), 460--475.

\bibitem{W}
S.~O.~Warnaar,
{\em Fermionic solution of the Andrews--Baxter--Forrester model. II.
Proof of Melzer's polynomial identities},
J. Stat. Phys. {\bf 84} (1996), 49--83.

\bibitem{Wb}
S.~O.~Warnaar,
{\em The Andrews--Gordon identities and $q$-multinomial coefficients},
q-alg/9601012. To appear in Commun. Math. Phys.


\end{thebibliography}
\end{document}